\documentclass[12pt]{article}

\makeatletter
\newfont{\footsc}{cmcsc10 at 8truept}
\newfont{\footbf}{cmbx10 at 8truept}
\newfont{\footrm}{cmr10 at 10truept}
\renewcommand{\ps@plain}{%
\renewcommand{\@oddfoot}{\footsc the electronic journal of combinatorics
  {\footbf 13} (2006), \#R1\hfil\footrm\thepage}}
\makeatother
\usepackage{latexsym}
\usepackage{amssymb}
\usepackage{amsmath}
\usepackage{fancybox}

\newcommand{\set}[1]{\left\{\,#1\,\right\}}

\newcommand{\botlabel}[2]{\mathrel{\mathop{\mathstrut{#1}}\limits
 _{\/
 \makebox[5pt]{\parbox[t]{80pt}{\centering
   $\uparrow$\par #2\par }}}}}
\DeclareMathOperator{\E}{E}

\title{Combinatorics of Partial Derivatives}

\author{Michael Hardy \\
\small School of Mathematics\\[-0.8ex]
\small University of Minnesota, Minneapolis, MN 55455, USA\\[-0.8ex]
\small \texttt{hardy@math.umn.edu}}

\date{\small 
Submitted: Oct 1, 2005;  Accepted: Dec 24, 2006; Published: Jan 7, 2006}

\begin{document}
\maketitle

\begin{abstract}
The natural forms of the Leibniz rule for the $k$th
derivative of a product and of Fa\`a di Bruno's formula
for the $k$th derivative of a composition involve the
differential operator $\partial^k/\partial x_1 \cdots \partial x_k$
rather than $d^k/dx^k$, with no assumptions about
whether the variables $x_1,\dots,x_k$ are all
distinct, or all identical, or partitioned into
several distinguishable classes of indistinguishable
variables.  Coefficients appearing in forms of these
identities in which some variables are indistinguishable
are just multiplicities of indistinguishable terms (in
particular, if all variables are distinct then all
coefficients are 1).  The computation of the multiplicities
in this generalization of Fa\`a di Bruno's formula is a
combinatorial enumeration problem that, although completely
elementary, seems to have been neglected.  We apply the
results to cumulants of probability distributions.
\end{abstract}

\section{Introduction}\label{intro}

Both the well-known Leibniz rule
\begin{equation}\label{Leibniz}
\frac{d^k}{dx^k}(uv)=\sum_{\ell=0}^k{k \choose \ell}
\frac{d^\ell u}{dx^\ell}\cdot\frac{d^{k-\ell} v}{dx^{k-\ell}}
\end{equation}
and the celebrated formula of Francesco Fa\`a di Bruno
\begin{equation}\label{Faa.di.Bruno}
{d^k \over dx^k} f(y)=\sum
{k! \over 1!^{m_1}\cdots k!^{m_k} m_1!\cdots m_k!}
f^{(m_1+\cdots+m_k)}(y)\prod_{j\,:\,m_j\neq 0}
{d^{m_j} y \over dx^{m_j}}
\end{equation}
(where the sum is over all $k$-tuples $(m_1,\dots,m_k)$
of non-negative integers satisfying the constraint
$m_1+2m_2+3m_3+\cdots+km_k=k$) are formulas for $k$th
derivatives of functions of functions of $x$.  That is
what the left sides of these identities share in common.
The right sides of both identities are sums whose terms
have products of higher derivatives with respect to $x$
as factors.

All mathematicians know the combinatorial interpretation
of the coefficients in the Leibniz rule (the number of
size-$\ell$ subsets of a size-$k$ set), and all
combinatorialists know the combinatorial interpretation
of the coefficients in Fa\`a di Bruno's formula (the number
of partitions of a size-$k$ set into $m_j$ parts of size
$j$, for $j=1,\dots,k$).  However, the following two
points appear not to be \mbox{widely known}:
\begin{enumerate}

\item The natural form of these identities involves
the differential operator
$$
\frac{\partial^k}{\partial x_1\,\cdots\,\partial x_k}
$$
instead of $d^k/dx^k$.  In that form, all coefficients
on the right sides are 1.

\item There should be no assumptions about whether the
variables $x_1,\dots,x_k$ are all distinct, or all identical,
or partitioned into several distinguishable classes of
indistinguishable variables.  When some variables become
indistinguishable, so do some of the terms on the right
sides of the identities.  Indistinguishable terms then
get collected, so that each constant coefficient of a
term on the right side is that term's multiplicity.
When {\em all} of the variables are indistinguishable,
then the multiplicities are the coefficients in
(\ref{Leibniz}) and (\ref{Faa.di.Bruno}) above.

\end{enumerate}

We will call the above Point 1 and Point 2.

Finding the multiplicities in Point 2 applied to (2) is
a combinatorial problem that may have escaped explicit
treatment until the present paper, in which the solution
is Proposition 4.  The problem is that of enumeration of
what we will call ``collapsing partitions''.

As an example of Point 2, if $x_2$ and $x_3$ collapse
into two indistinguishable variables called $x_2$, so that
$\partial^3/(\partial x_1\,\partial x_2\,\partial x_3)$
becomes $\partial^3/(\partial x_1\,\partial x_2^2)$,
then the two-term sum
$$
{\partial y \over \partial x_2}
\cdot
{\partial^2 y \over \partial x_1\, \partial x_3}
+
{\partial y \over \partial x_3}
\cdot
{\partial^2 y \over \partial x_1\, \partial x_2}
$$
collapses to the term
$$
2\cdot {\partial y \over \partial x_2}
\cdot
{\partial^2 y \over \partial x_1\, \partial x_2}
$$
with multiplicity 2.

The chain rule and the product rule are enough to entail
that the coefficients must be positive integers.  But
without Point 1 and Point 2, it is not obvious what, if
anything, they enumerate.

Two papers, Constantine and Savits \cite{Constantine+Savits}
and Leipnik and Reid \cite{Leipnik+Reid}, give an identity expressing 
$(\partial^{k_1+\cdots+k_n}/\partial x_1^{k_1}
\cdots\partial x_n^{k_n})f(y)$
as a linear combination of products of derivatives
of $f(y)$ with respect to $y$ and of $y$ with respect
to the independent variables.
But neither of those sources mentions that as more and
more variables become indistinguishable the identity does
not change except in the collection of newly
indistinguishable terms.  Without that observation, the
combinatorial content of the problem is invisible.
Leipnik and Reid in \cite{Leipnik+Reid}, p.~1, wrote,
``Obviously, `pure' derivatives, such as
$\frac{\partial^4 G(z_1,z_2)}{\partial z_1^4}$
are easier to deal with than mixed derivatives like
$\frac{\partial^4 G}{\partial z_1^2\,\partial z_2^2}$.''
But from our point of view, it will be {\em maximally}
``mixed'' derivatives like
$\partial^4 G/(\partial z_1\,\partial z_2\, \partial z_3\,\partial z_4)$
that are the easiest and most basic.

Proposition 2 of this paper is partially anticipated by
Terry Speed in \cite{Speed}, page 382.  That paper gives
only the special case in which $f$ is the exponential
function, in which the derivatives appear as coefficients
of power series, and is stated in an inconspicuous
and somewhat tangential way that mixes it so thoroughly
with the theory of cumulants in probability theory that
it can be understood only by understanding what the paper
is saying about cumulants.  Speed wrote:
``\dots the general results are most transparent when all
\dots variables under discussion are taken to be distinct.''
That remark played a role in inspiring this paper.
Its influence will be seen not only in our Proposition 1,
but also in our treatment of product rules--a topic not
directly relevant to that of Speed's paper and not mentioned
there.  Speed went on:
``The identification of some or all [variables] at a later
stage merely introduces extra factors, and at times these
multiplicities are not particularly easy to calculate.''
The multiplicities are given by our Proposition 4.

Unlike \cite{Constantine+Savits} and \cite{Leipnik+Reid}, 
Warren Johnson \cite{Johnson} states a version of Fa\`a
di Bruno's formula that is explicit about the combinatorial
meaning of the coefficients.  But Johnson treats only
functions of one variable and gives nothing like
Proposition 2 of the present paper.  The same is true
of the ``compositional formula'' on page 3 of
Richard P.~Stanley's treatise \cite{Stanley}.
Like Speed, Stanley gives a power-series version of the
formula.  He mentions the name of Fa\`a di Bruno only
in endnotes.  One also finds a very combinatorics-flavored
views of Fa\`a di Bruno's formula in \cite{Zeilberger}.
Other variations on the theme are in \cite{Mishkov},
\cite{Noschese+Ricci}, and \cite{Yang}.

We conclude this paper with the application of
the results to cumulants.

\section{Partial derivatives and partitions of sets}

In the identity
\begin{equation}\label{firstidentity}
\begin{array}{rcl}
\dfrac{\partial^3}{\partial x_1\, \partial x_2\, \partial x_3}e^y
& = & e^y\left(\dfrac{\partial^3 y}{
\partial x_1\, \partial x_2\, \partial x_3}
\,+\,\dfrac{\partial y}{\partial x_1}
\cdot\dfrac{\partial^2 y}{\partial x_2\, \partial x_3}
+\,\dfrac{\partial y}{\partial x_2}
\cdot\dfrac{\partial^2 y}{\partial x_1\, \partial x_3}\right. \\  \\
& & \qquad\left.+\,\dfrac{\partial y}{\partial x_3}
\cdot\dfrac{\partial^2 y}{\partial x_1\, \partial x_2}
\,+\,\dfrac{\partial y}{\partial x_1}
\cdot\dfrac{\partial y}{\partial x_2}
\cdot\dfrac{\partial y}{\partial x_3} \right),
\end{array}
\end{equation}
\enlargethispage{5mm}
where $y$ is a function of $x_1, x_2, x_3$,
the terms correspond in an obvious way to the
five partitions of the set $\{\,1,2,3\,\}$.
We will see that this holds generally:
the partial derivative
$(\partial^n/\partial x_1\cdots \partial x_n)e^y$
is $e^y$ times the sum whose terms correspond in
just this way to the partitions of the set $\set{1,\dots,n}$.
Using the notation (for example)
$\partial^3 y \left/ \prod_{j\in\set{2,4,9}}\partial x_j\right.$
to mean
$\left.\partial^3 y \right/ \partial x_2\,\partial x_4\,\partial x_9$,
we can say
\begin{equation}\label{exponentialformula}
{\partial^n \over \partial x_1\cdots\partial x_n}e^y
=e^y\sum_{\pi} \prod_{B\in\pi}
\dfrac{\partial^{\left|B\right|} y}{
\prod_{j\in B}\partial x_j}
\end{equation}
where the sum is over all partitions $\pi$ of the set
$\set{1,\dots,n}$ and the product is over all of the
parts $B$, or ``blocks'' as we will call them, of the
partition $\pi$, and we denote the number of members
of any set $S$  by $\left|S\right|$.

If we have $f(y)$ instead of $e^y$, then the
orders of the derivatives of $f$ must be mentioned.
The order of each derivative of $f$ is just the number
of blocks in the partition.  This is the first result
that we will prove (in Section \ref{Proofs}):
                                                                                
{\bf Proposition 1.}\label{1}
\begin{equation}\label{main}
{\partial^n \over \partial x_1\cdots\partial x_n}f(y)
=\sum_{\pi} f^{(\left|\pi\right|)}(y) \prod_{B\in\pi}
\dfrac{\partial^{\left|B\right|} y}{
\prod_{j\in B}\partial x_j}.
\end{equation}

{\bf Example 1.}
\begin{eqnarray*}
& & {\partial^3 \over \partial x_1\, \partial x_2\, \partial x_3}f(y)
=\underbrace{f'(y){\partial^3 y
\over \partial x_1\, \partial x_2\, \partial x_3}}_{
\text{1 block; 1st derivative of }f.} \\  \\
& + & \,\underbrace{f''(y) \left( {\partial y \over \partial x_1}
\cdot{\partial^2 y \over \partial x_2\, \partial x_3}
\,+\,{\partial y \over \partial x_2}
\cdot{\partial^2 y \over \partial x_1\, \partial x_3}
+\, {\partial y \over \partial x_3}
\cdot{\partial^2 y \over \partial x_1\, \partial x_2}\right)}_{
\text{2 blocks in each partition; 2nd derivative of }f.} \\  \\
& + & \,\underbrace{f'''(y) {\partial y \over \partial x_1}
\cdot{\partial y \over \partial x_2}
\cdot{\partial y \over \partial x_3}}_{
\text{3 blocks; 3rd derivative of }f.}.
\end{eqnarray*}

{\bf Proposition 2.}\label{2}
If some of the $x_i$s become indistinguishable,
then so do the corresponding terms in the sum;
nothing else changes.

This will also be proved in Section \ref{Proofs}.

{\bf Example 2.} Suppose that in Example 1, the two
variables $x_2$ and $x_3$ become indistinguishable from
each other.  Call them both $x_2$.  Then we have
\begin{eqnarray*}
{\partial^3 \over \partial x_1\, \partial x_2^2}f(y)
& = & \overbrace{f'(y){\partial^3 y
\over \partial x_1\, \partial x_2^2}}^{
\text{1 block; 1st derivative of }f.} \\  \\
& & +\,\overbrace{f''(y) \left( {\partial y \over \partial x_1}
\cdot{\partial^2 y \over \partial x_2^2}\right.
\,+\,2\cdot  \botlabel{
{\underbrace{{\partial y \over \partial x_2}
\cdot{\partial^2 y \over \partial x_1\, \partial x_2}}}
}{\ovalbox{$\begin{array}{c}
\text{The multiplicity} \\
\text{of this term is }2.
\end{array}$}} \left.\phantom{\partial^1 \over \partial_1}
\hspace{-0.5cm}
\right) }^{\text{
2 blocks in each partition; 2nd derivative of }f.} \\  \\
& & +\,\overbrace{f'''(y) {\partial y \over \partial x_1}
\cdot\left({\partial y \over \partial x_2}\right)^2}^{
\text{3 blocks; 3rd derivative of }f.}.
\end{eqnarray*}

The multiplicity mentioned above is how many formerly
distinguishable terms get collected to form that term.
The problem of finding such multiplicities is treated
in Section \ref{multisets}.  Applying the language of
that section to Example 2, we would ask: how many partitions
of the set $\set{1,2,3}$ collapse to the partition
$\set{2}+\set{1,2}$ of the multiset $\set{1,2,2}$ when the
set $\set{1,2,3}$ collapses to the multiset $\set{1,2,2}$,
i.e., when the members 2 and 3 become indistinguishable?
The answer in this case is 2.

\section{Proofs of the first two Propositions}\label{Proofs}

{\bf Proof of Proposition 1.}
This proof relies on this simple standard algorithm
for converting a list of all partitions of
$\set{1,\dots,n}$ into a list of all partitions
of $\set{1,\dots,n+1}$:

\begin{enumerate}\label{algorithm}

\item To each partition of $\set{1,\dots,n}$, add the
1-member-set $\set{n+1}$ as a new block.  This gives
a list of some of the partitions of $\set{1,\dots,n+1}$.

\item To each block of each partition of $\set{1,\dots,n}$,
add $n+1$ as a new member of the block.  This gives a list
of $\sum_\pi \left|\pi\right|$ additional partitions of
$\set{1,\dots,n+1}$.

\end{enumerate}

The union of these two lists clearly contains all
partitions of $\set{1,\dots,n+1}$.

In particular, it works when $n=0$, since the empty set
has exactly one partition\footnote{Perhaps as a result
of studying set theory, I was surprised when I learned
that some respectable combinatorialists consider such
things as this to be mere convention.  One of them even
said a case could be made for setting the number of
partitions to 0 when $n=0$.  By stark contrast, Gian-Carlo Rota
wrote in \cite{Rota2}, p.~15, that ``the kind of mathematical
reasoning that physicists find unbearably pedantic''
leads not only to the conclusion that the elementary
symmetric function in no variables is 1, but straight
from there to the theory of the Euler characteristic,
so that ``such reasoning does pay off.''  The only other
really sexy example I know is from applied statistics:
the non-central chi-square distribution with zero degrees
of freedom, unlike its ``central'' counterpart, is non-trivial.}.
That establishes the basis for a proof by mathematical
induction on $n$.

Next we use the Proposition in case $n$ to prove the
Proposition in case $n+1$.
\begin{eqnarray*}
& & {\partial^{n+1} \over \partial x_1\cdots\partial x_{n+1}}f(y)
\\  \\
& = & \sum_{\pi}
{\partial \over \partial x_{n+1}}\left[f^{(\left|\pi\right|)}(y)
\prod_{B\in\pi}
\dfrac{\partial^{\left|B\right|} y}{
\prod_{j\in B}\partial x_j}\right] \\  \\
& = & \sum_{\pi}\left[
f^{(\left|\pi\right|+1)}(y)
{\partial y\over \partial x_{n+1}}
\prod_{B\in\pi}
\dfrac{\partial^{\left|B\right|} y}{
\prod_{j\in B}\partial x_j}
+f^{(\left|\pi\right|)}(y)
{\partial \over \partial x_{n+1}} \prod_{B\in\pi}
\dfrac{\partial^{\left|B\right|} y}{
\prod_{j\in B}\partial x_j}\right] \\  \\
& = & \sum_{\pi}\left[
f^{(\left|\pi\right|+1)}(y)
{\partial y\over \partial x_{n+1}}
\prod_{B\in\pi}
\dfrac{\partial^{\left|B\right|} y}{
\prod_{j\in B}\partial x_j}   \right. \\  \\
& & \qquad\qquad                   \left.
+f^{(\left|\pi\right|)}(y)
\sum_{B\in\pi}
\left(\dfrac{\partial^{\left|B\right|+1}y}{
\partial x_{n+1}\prod_{j\in B} \partial x_j}
\cdot\prod_{C\in\pi\,:\,C\neq B}
\dfrac{\partial^{\left|C\right|}y}{\prod_{j\in C}\partial x_j}
\right)\right].
\end{eqnarray*}
Inside the last square brackets is a sum of two terms.
The first term corresponds to step 1 in our algorithm:
we have added $\set{n+1}$ as a new block to our partition
$\pi$ of $\set{1,\dots,n}$, getting a partition with
$\left|\pi\right|+1$ blocks, of the set $\set{1,\dots,n+1}$.
The second corresponds to step 2 in our algorithm:
we have added $n+1$ to each block of our partition $\pi$
of $\set{1,\dots,n}$, getting a partition with
$\left|\pi\right|$ blocks, of the set $\set{1,\dots,n+1}$.
We now have a sum over all partitions of $\set{1,\dots,n+1}$,
each partition being represented as a product of partial
derivatives, each partial derivative representing a block.
And for each partition there is a factor $f^{(\bullet)}(y)$,
the order of the derivative being the number of blocks
of the partition.  This proves case $n+1$, and the proof
by induction on $n$ is complete. $\blacksquare$

{\bf Proof of Proposition 2.} Observe that if, in the
argument above, we had differentiated at the $(n+1)$th
step with respect to $x_k$ for some $k\in\set{1,\dots,n}$,
rather than with respect to $x_{n+1}$, then nothing would
change except that some formerly distinguishable terms
would become indistinguishable. $\blacksquare$

\section{Multisets and collapsing partitions}\label{multisets}

\subsection{Definitions and conventions}

The first two bullet points and the fifth in the
definition below are standard but make clear which
notational conventions we will follow.  The third and
fourth may be less standard.
\begin{itemize}

\item A {\bf multiset} is a ``set with multiplicities'',
i.e., positive integers, assigned to each member $x$,
thought of as the number of times $x$ occurs as a member.
We will write $\left\{\right.\underbrace{x_1,\dots,x_1}_{m_1},
\dots,\underbrace{x_n,\dots,x_n}_{m_n}\left.\right\}$,
indicating multiplicities with underbraces, or, in simple
cases, for example $\set{a,a,a,b,b,c,c,c,c,c}$, the multiplicity
being the number of times the member is named.  In
particular, we identify every set with a multiset in
which every multiplicity is 1.

Perhaps the most widely known example is the multiset of
prime factors of a natural number: each prime factor
has a multiplicity.

\item The {\bf size} $\left|S\right|$ of a multiset
$S$ is the sum of the multiplicities.

\item The {\bf sum} of multisets is given by term-by-term
addition of multiplicities:
\enlargethispage{5mm}
\begin{eqnarray*}
& & \left\{\right.\underbrace{x_1,\dots,x_1}_{\ell_1},
\dots,\underbrace{x_n,\dots,x_n}_{\ell_n}\left.\right\}
+ \left\{\right. \underbrace{x_1,\dots,x_1}_{m_1},
\dots,\underbrace{x_n,\dots,x_n}_{m_n} \left.\right\} \\  \\
& = & \left\{\right.\underbrace{x_1,\dots,x_1}_{\ell_1+m_1},
\dots,\underbrace{x_n,\dots,x_n}_{\ell_n+m_n} \left.
\right\}.
\end{eqnarray*}
Only when sets are disjoint is their sum the same
as their union.

\item
A {\bf partition} of a multiset expresses that
multiset as a sum of multisets.

\item A {\bf partition} of a positive integer
expresses that integer as a sum of positive integers.

\end{itemize}
The next proposition is trivial but crucial.

{\bf Proposition 3.}

\begin{itemize}

\item The concept of {\em partition of a set}
is a special case of that of {\em partition of
a multiset}.

\item If we identify any multiset in which
``all members are equal'' (i.e., there is just one member,
whose multiplicity may be any positive integer) with the
multiplicity of that one member (for example, the multiset
$\set{a,a,a}$ is identified with the number 3), then the
concept of {\em partition of an integer} becomes a special
case of that of {\em partition of a multiset}.

\end{itemize}

\subsection{Collapsing partitions}

If the members 1, 2, 3, 4 of the set $\set{1,2,3,4,5,6,7,8}$
are made indistinguishable from each other and are called ``1'',
and 5 and 6 are made indistinguishable from each other and are
called ``5'', then we say that the set $\set{1,2,3,4,5,6,7,8}$
has ``collapsed'' to the multiset $\set{1,1,1,1,5,5,7,8}$.
Then we can ask: how many set-partitions collapse to the
multiset-partition
$$
\set{1,1,5}+\set{1,1,5}+\set{7,8}\text{?}
$$
It is a simple exercise to find that the answer is 6.
Consequently, via the correspondence
\begin{equation}\label{correspondence}
\tau=\left\{\right.\underbrace{1,\dots,1}_{k_1},
\dots,\underbrace{n,\dots,n}_{k_n}\left.\right\}
\longleftrightarrow
{\partial^{k_1+\cdots+k_n} \over \partial x_1^{k_1}
\cdots \partial x_n^{k_n}}=\partial_\tau
\end{equation}
between multisets and partial differential operators,
Proposition 2 entails that the expansion of the partial
derivative
\begin{equation}\label{unexpanded.derivative}
\partial_{\set{1,1,1,1,5,5,7,8}}f(y)
={\partial^8 \over \partial x_1^4\,\partial x_5^2\,
\partial x_7\,\partial x_8}f(y)
\end{equation}
contains (among many others) this term:
\begin{equation}\label{particular.term}
6f'''(y)\left(\partial_{\set{1,1,5}} y\right)^2\,\partial_{\set{7,8}} y
=6f'''(y)\left({\partial^3 y
\over \partial x_1^2\,\partial x_5}\right)^2
\cdot {\partial^2 y \over \partial x_7\,\partial x_8}
\end{equation}
(the order of the derivative of $f$ is 3 because that is
how many blocks are in this partition).

An extreme case of ``collapsing'' is exemplified by the
question: How many partitions of the set $\set{1,2,3,4,5,6,7,8}$
collapse to the partition $3+3+2$ of the number 8 when all 8
members of the set $\set{1,2,3,4,5,6,7,8}$ collapse into
indistinguishability?  Again, a simple exercise shows that
the answer is 280.  In this extreme case where all members
of the set become indistinguishable and partitions of the
multiset become partitions of an integer, the answer to the
problem of enumeration of collapsing partitions is well known
to be given by the coefficients in Fa\`a di Bruno's formula
-- in this case by the coefficient of
$$
f'''(y) \left({d^3 y \over dx^3}\right)^2 {d^2 y \over dx^2}
$$
in the expansion of $(d^8 / dx^8) f(y)$.

In general, we have this result:

{\bf Corollary to Propositions 1 and 2.}
Let $\tau=\left\{\right.\underbrace{1,\dots,1}_{k_1},
\dots,\underbrace{n,\dots,n}_{k_n}\left.\right\}$.  Use the
notation introduced in (\ref{correspondence}) above.
Then
$$
\partial_\tau f(y)
= {\partial^{k_1+\cdots+k_n} \over \partial x_1^{k_1}
\cdots \partial x_n^{k_n}}f(y)
=\sum_{\tau_1+\tau_2+\cdots}
M f^{(\bullet)}(y)\cdot \left(\partial_{\tau_1} y
\cdot \partial_{\tau_2} y\cdot\,\cdots\right),
$$
where the sum is over all partitions $\tau_1+\tau_2+\cdots$
of the multiset $\tau$, the order of the derivative
$f^{(\bullet)}(y)$ is the number of terms in the partition
$\tau_1+\tau_2+\cdots$, and the multiplicity $M$ is the number
of partitions of the set $\set{1,2,3,\dots,k_1+\cdots+k_n}$
that collapse to the partition $\tau_1+\tau_2+\cdots$ of the
multiset $\tau$ when the set $\set{1,2,3,\dots,k_1+\cdots+k_n}$
collapses to the \mbox{multiset $\tau$}.

In order to use it in the next result, we
introduce a convention:

{\bf Notational convention.}  For any multiset $\sigma$
let $\sigma!!$ denote the product of the factorials of
the multiplicities of the members of $\sigma$.
For example, $\set{1,1,1,1,2,2,2}!!=4!3!=144$.

The next result will be proved in Section \ref{multiplicitiesproof}:

{\bf Proposition 4.}

Let $\tau=\left\{\right. \underbrace{1,\dots,1}_{k_1},
\dots,\underbrace{n,\dots,n}_{k_n}\left.\right\}$.
Consider a partition
$$
\tau=\underbrace{\tau_1+\cdots+\tau_1}_{m_1}
+\underbrace{\tau_2+\cdots+\tau_2}_{m_2}
+\cdots
$$
in which $\tau_1$, $\tau_2$, \dots are {\em all distinct}
(so that $m_1, m_2, m_3,\dots$ are multiplicities of
yet another sort).  Denote this by
$$
\tau=m_1\tau_1+m_2\tau_2+\cdots.
$$
Then the number of partitions of the set
$\set{1,2,3,\dots,k_1+\cdots+k_n}$ that collapse to
the partition $m_1\tau_1+m_2\tau_2+m_3\tau_3+\cdots$
of the multiset $\tau$ when the set
$\set{1,2,3,\dots,k_1+\cdots+k_n}$ collapses to
the multiset $\tau$ is
\begin{equation}\label{multiplicities}
\frac{k_1!\cdots k_n!}{
\tau_1!!^{m_1}
\tau_2!!^{m_2}
\tau_3!!^{m_3}\cdots
m_1! m_2! m_3! \cdots }.
\end{equation}

\subsection{The most extreme case}

In the most extreme case of indistinguishability of independent
variables, the operator $\partial^k / \partial x_1 \cdots \partial x_k$
collapses to $d^k/dx^k$:
\begin{equation}\label{collapsed}
{d^k \over dx^k}f(y)=\sum_{\pi} f^{(\left|\pi\right|)}(y) \prod_{B\in\pi}
\left(d \over dx\right)^{\left|B\right|}y,
\end{equation}
where, again, the sum is over all partitions $\pi$
of $\set{1,\dots,k}$.
In this sum, two terms are indistinguishable whenever
two partitions of the set $\{\,1,\dots,k\,\}$ both collapse
to the same partition of the integer $k$ when {\em all} of
the members of $\{\,1,\dots,k\,\}$ become indistinguishable.

In this extreme case, the multiplicities are given
by the classic formula of Francesco
Fa\`a di Bruno (\ref{Faa.di.Bruno}).
Francesco Fa\`a di Bruno (1825 -- 1888) was
(in chronological order) a military officer,
a mathematician, and a priest.  He published
this formula in \cite{Bruno1} and \cite{Bruno2} and
was posthumously beatified by the Pope\footnote{An
anonymous referee suggests that that pontiff may have
been influenced more by Fa\`a di Bruno's charitable
than mathematical work.}.

Alex Craik's {\em Prehistory of Faa di Bruno's Formula}
\cite{Craik} points out that Fa\`a di Bruno was
anticipated in 1800 by L.F.A.~Arbogast; see \cite{Arbogast}.

Although (\ref{Faa.di.Bruno})
is the well-known form of this identity, Warren P.~Johnson
has written in \cite{Johnson} (bottom of page 231), that
(\ref{collapsed}) ``is really the fundamental form''
of Fa\`a di Bruno's formula.  We propose that the conjunction
of our first two Propositions is more fundamental.

\subsection{Conservation of Bell numbers}
                                                                                
By now it should be clear that, when the derivative
$$
{\partial^n\over \cdots\,\cdots}f(y)
$$
is expanded as a sum in terms of derivatives of $f(y)$
with respect to $y$ and derivatives of $y$ with respect
to whichever independent variables appear, then the sum of
the coefficients is always the number $B_n$ of partitions
of a set of $n$ members.  This is called the $n$th
Bell number, in honor of Eric Temple Bell.  The first
several of these are $B_0=1, B_1=1, B_2=2, B_3=5,
B_4=15, B_5=52, B_6=203, B_7=877, B_8=4140,\dots$.
For an account of these numbers, see \cite{Rota}.
As more and more of the $n$ independent variables get
identified with each other, the sum of the multiplicities
never changes.

\section{Proof of Proposition 4}\label{multiplicitiesproof}

Imagine $k_1+\cdots+k_n$ Scrabble tiles.
On the first $k_1$ of these, the number ``1'' is written;
on the next $k_2$ of them, ``2'' appears; and so on.
These are partitioned into $m_1$ copies of the multiset
$\tau_1$, $m_2$ copies of $\tau_2$, and so on.

Permuting the $k_1$ ``1''s or the $k_2$ ``2''s, etc., 
does not alter the set of $m_1+\cdots+m_n$ ``words''.
Thus there are $k_1!\cdots k_n!$ permutations of the
$k_1+\cdots+k_n$ tiles representing the partition
$$
\tau=m_1\tau_1+m_2\tau_2+\cdots.
$$
Permuting the identical elements within any block of this
partition does not alter which partition of $\tau$ we have,
and therefore the product $k_1!\cdots k_n!$ gets divided by
$\tau_1!!^{m_1} \tau_2!!^{m_2}\cdots$.
Neither does permuting the $m_i$ blocks identical to $\tau_i$,
and therefore it also gets divided by $m_1!\cdots m_n!$.
$\blacksquare$

\section{Product rules}

In the identity
\begin{eqnarray*}
& & {\partial^3 \over \partial x_1\,\partial x_2\,\partial x_3}(uv) \\  \\
& = & u\cdot{\partial^3 v \over \partial x_1\,\partial x_2\,\partial x_3}
+{\partial u \over \partial x_1}
\cdot {\partial^2 v \over \partial x_2\,\partial x_3}
+{\partial u \over \partial x_2}
\cdot {\partial^2 v \over \partial x_1\,\partial x_3}
+{\partial u \over \partial x_3}
\cdot {\partial^2 v \over \partial x_1\,\partial x_2} \\  \\
& & +{\partial^2 u \over \partial x_1\,\partial x_2}
\cdot{\partial v \over \partial x_3}
+{\partial^2 u \over \partial x_1\,\partial x_3}
\cdot{\partial v \over \partial x_2}
+{\partial^2 u \over \partial x_2\,\partial x_3}
\cdot{\partial v \over \partial x_1}
+{\partial^3 u \over \partial x_1\,\partial x_2\,\partial x_3}\cdot v,
\end{eqnarray*}
the terms correspond in an obvious way to the eight
subsets of the set $\set{1,2,3}$.  This exemplifies
the first part of our next result.

{\bf Proposition 5.}
\begin{itemize}
\item
$$
{\partial^n \over \partial x_1 \cdots \partial x_n}(uv)
=\sum_{S}\dfrac{\partial^{(\left|S\right|)} u}{
\prod_{j\in S} \partial x_j}
\cdot \dfrac{\partial^{(n-\left|S\right|)}v}{
\prod_{j\not\in S} \partial x_j},
$$
where the index $S$ runs through the set of all
subsets of $\set{1,\dots,n}$.

\item If some of the variables become indistinguishable,
then so do some of the terms in the sum; nothing else
changes.

\end{itemize}

{\bf Example 4.}
\begin{eqnarray*}
{\partial^3 \over \partial x_1\,\partial x_2^2}(uv)
& = & u\cdot{\partial^3 v \over \partial x_1\,\partial x_2^2}
+{\partial u \over \partial x_1}
\cdot {\partial^2 v \over \partial x_2^2}
+2\cdot{\partial u \over \partial x_2}
\cdot {\partial^2 v \over \partial x_1\,\partial x_2} \\  \\
& & +2\cdot {\partial^2 u \over \partial x_1\,\partial x_2}
\cdot{\partial v \over \partial x_2}
+{\partial^2 u \over \partial x_2^2}
\cdot{\partial v \over \partial x_1}
+{\partial^3 u \over \partial x_1\,\partial x_2^2}\cdot v.
\end{eqnarray*}

In this case the solution of the combinatorial
problem is simpler:

{\bf Proposition 6.}
$$
{\partial^{k_1+\cdots+k_n} \over \partial x_1^{k_1}\,\cdots\,
\partial x_n^{k_n}} (uv) 
=\sum_{\ell_1=0}^{k_1}\cdots\sum_{\ell_n=0}^{k_n}
{k_1 \choose \ell_1}\cdots {k_n \choose \ell_n}
{\partial^{\ell_1+\cdots+\ell_n} u
\over \partial x_1^{\ell_1}\,\cdots\,\partial x_n^{\ell_n}}
\cdot {\partial^{k_1-\ell_1+\cdots+k_n-\ell_n} v
\over \partial x_1^{k_1-\ell_1}\,\cdots
\,\partial x_n^{k_n-\ell_n}}.
$$
In the most extreme case, all of the independent
variables become indistinguishable, and we have
the familiar Leibniz rule (\ref{Leibniz}).
The proofs of Propositions 5 and 6 are left as exercises.

As far as the present writer knows, the second part of
Proposition 5 is new.  It identifies the easy combinatorial
problem that Proposition 6 solves.
Proposition 6 can be found in both \cite{Constantine+Savits}
and \cite{Comtet}, p.~131, but without the combinatorial
interpretation of the coefficients.  The first part of
Proposition 5 is an important special case of Proposition 6,
and we do not know of any earlier explicit statement of it
than that in the present paper (a referee says ``it could
be just about anywhere'' and I have not succeeded in proving
otherwise).

\section{Cumulants}
The omitted fragment represented by the second ellipsis ``\dots''
in the first quote from Terry Speed in Section \ref{intro}
is the word {\em random}.  That is because Speed's topic
is that of cumulants of random variables.

For positive integers $n$, the $n$th cumulant functional
$\kappa_n$ assigns a real number $\kappa_n(X)$ to real-valued
random variables $X$.  Let $\mu=\E(X)$ be the expected value
of $X$. Then the $n$th central moment of $X$ is $\E((X-\mu)^n)$.
For $n\geq 2$, the $n$th cumulant shares with the $n$th central
moment the properties of $n$th-degree homogeneity and
translation-invariance:
\begin{eqnarray*}
\kappa_n(cX) & = & c^n \kappa_n(X), \\
\kappa_n(X+c) & = & \kappa_n(X).
\end{eqnarray*}
Moreover, if the random variables $X_1,\dots,X_m$
are independent, then
$$
\kappa_n(X_1+\cdots+X_m)
=\kappa_n(X_1)+\cdots+\kappa_n(X_m).
$$
The $n$th central moment has this additivity property
{\em only} when $n\leq 3$.  In fact, when $n=$ either 2
or 3, the cumulant is just the central moment.
When $n=1$, then the cumulant is the expected value.
For all $n$, the $n$th cumulant is an $n$th-degree
polynomial in the first $n$ moments.

Cumulants were introduced in the 19th century in \cite{Thiele}
by the Danish actuary Thorvald Thiele, who called them
half-invariants; an English translation was published in
1931; see \cite{Thiele2}.  They were first publicly given
the name {\em cumulants} in 1931 by the statisticians
Ronald Fisher and John Wishart in \cite{Fisher+Wishart},
the name having been suggested to Fisher in private
correspondence from the statistician Harold Hotelling.
There are also {\em joint cumulants}
$$
\kappa(X_1,\dots,X_n).
$$
When $n=2$, this is just the covariance.
All probabilists and statisticians know that the covariance
between a random variable and itself is its variance.  A
similar thing happens with joint cumulants: When the $n$
random variables collapse into indistinguishability, then
the joint cumulant coincides with the $n$th cumulant of
one random variable:
$$
\kappa(\,\underbrace{X,\dots,X}_{n}\,)=\kappa_n(X).
$$
Beyond this talk of ``collapsing into indistinguishability'',
a parallel between cumulants and the partial derivatives
treated in the foregoing sections is seen in the identity
that expresses the $n$th raw moment $\E(X^n)$ (not the $n$th
{\em central} moment) in terms of the first $n$ cumulants:
\begin{equation}\label{moments.from.cumulants}
\E(X^n)=\sum_\pi \prod_{B\in\pi}\kappa_{\left|B\right|}(X).
\end{equation}
For random variables $X_1,\dots,X_n$, a similar identity holds:
\begin{equation}\label{joint.moments.from.joint.cumulants}
\E(X_1\cdots X_n)=\sum_\pi\prod_{B\in\pi} \kappa(X_i : i \in B).
\end{equation}
The identities (\ref{moments.from.cumulants}) and
(\ref{joint.moments.from.joint.cumulants}) completely
characterize all of the cumulant functionals.

That is how Terry Speed came to consider the question
of these multiplicities, of which all he wrote was that
they are not always easy to calculate.  Workers with
cumulants know what to do in the two opposite extreme
cases: when {\em none} of the random variables are
identified then all coefficients are equal to 1, and
when {\em all} are identified then the classic
Fa\`a di Bruno formula (\ref{Faa.di.Bruno}) gives the
coefficients.  But if one were to judge by Speed's comments,
they might seem at something of a loss in cases intermediate
between ``all'' and ``none''.  That case is handled by our
Proposition 4.

Speed's topic and that of partial derivatives are not
two disparate applications of our Proposition 4.  The
joint cumulant may be characterized as the coefficient
of $t_1\cdots t_n$ in the power-series expansion of
$$
\log \E\left(\exp(t_1 X_1+\cdots+t_n X_n)\right)
$$
(at least in the case in which all moments exist).
More tersely stated, the joint cumulant-generating function
is the logarithm of the joint moment-generating function.
Since the coefficients of a power series of the form
$\sum_{i_1,\dots,i_n}c_{i_1,\dots,i_n}
t_1^{i_1}\cdots t_n^{i_n}/(i_1!\cdots i_n!)$ are its partial
derivatives at $t_1=\cdots=t_n=0$, Speed's topic is a special
case of ours.

\section{Acknowledgments}

Ezra Miller and Jay Goldman offered helpful suggestions.
Unlike some copyeditors, Ellen Stuttle found the actual
content interesting.  Alex Craik pointed out some
relevant references and reminded me of Arbogast's
anticipation of Fa\`a di Bruno's formula.  In addition
to some encouraging words, Doron Zeilberger pointed out
his paper \cite{Zeilberger}; which is even more ``discretely''
oriented than the present paper (which, the reader will have
noticed, disdains even to mention such ``analytical'' things
as differentiability).

\end{document}